\documentclass[aop,MSNbibl,seceqn,citesort,dvips]{arximspdf}
\usepackage{stmaryrd}
\usepackage{wasysym}
\usepackage{graphicx}

\doi{10.1214/11-AOP740} %kopijuoti is PTS
\volume{41}
\issue{5}
\pubyear{2013}
\firstpage{3261}
\lastpage{3283}

\makeatletter

\newcommand{\ext}{\operatorname{ext}}
\newcommand{\sT}{{\mathcal T}}
\newcommand{\oo}{\infty}
\newcommand{\De}{\Delta}
\newcommand{\rad}{\operatorname{rad}}
\newcommand{\si}{\sigma}
\newcommand{\pd}{\partial}
\newcommand{\es}{\varnothing}
\newcommand{\Om}{\Omega}
\newcommand{\Ga}{\Gamma}
\newcommand{\La}{\Lambda}
\newcommand{\sM}{{\mathcal M}}
\newcommand{\pc}{p_{\mathrm{c}}}
\newcommand{\lra}{\leftrightarrow}
\newcommand{\xlra}{\xleftrightarrow}
\newcommand{\nxlra}[1]{\mathrel{\mathpalette\concel{\xleftrightarrow{#1}}}}
\newcommand{\concel}[2]{\ooalign{$\hfil #1\mkern0mu/\hfil$\crcr$#1#2$}}

\newcommand{\xleftrightarrow}[1]{\mathrel{\vbox{\m@th\ialign{##\crcr
  $\hspace*{3pt}\hfil\scriptstyle\ #1\ \ \hfil$\crcr\noalign{\kern0.5pt\nointerlineskip}%
  \leftarrowfill$\hspace*{-2pt}\rightarrow$\crcr}}}}

\newcommand{\bp}{\mathbf{p}}
\newcommand{\squ}{{\square}}
\newcommand{\tri}{{\triangle}}
\newcommand{\hex}{{\hexagon}}
\newcommand{\Cv}{C_{\mathrm{v}}}
\newcommand{\Ch}{C_{\mathrm{h}}}
\newcommand{\eps}{\varepsilon}
\newcommand{\bpc}{{\bp_{\mathrm{c}}}}

\newcommand{\RR}{\mathbb{R}} % real numbers
\newcommand{\NN}{\mathbb{N}} % non-negative integers
\newcommand{\ZZ}{\mathbb{Z}} % integers
\newcommand{\PP}{\mathbb{P}} % probability
\newcommand{\HH}{\mathbb{H}}
\newcommand{\TT}{\mathbb{T}}
\newcommand{\LL}{\mathbb{L}}
\newcommand{\comp}{\circ}

\newcommand{\Su}{S^{\Yup}}
\newcommand{\Sd}{S^\Ydown}
\newcommand{\Tu}{T^\vartriangle}
\newcommand{\Td}{T^\triangledown}
\newcommand{\Ann}{\mathcal{A}}
\newcommand{\Lattice}{\mathbb{L}}
\newcommand{\Lat}{\Lattice}

\newtheorem{theorem}{Theorem}[section]
\newtheorem{lemma}[theorem]{Lemma}
\newtheorem{prop}[theorem]{Proposition}
\newtheorem{cor}[theorem]{Corollary}

\newproclaim{conj}[theorem]{Conjecture}

\newcommand{\Arm}{A}
\newcommand{\oArm}{\bar\Arm}
\newcommand{\oAnn}{\bar\Ann}
\newcommand{\eqv}{\asymp}
\newcommand{\logeqv}{\approx}

\makeatother

\begin{document}
\begin{frontmatter}

\title{Universality for bond percolation in~two~dimensions}
\runtitle{Universality for bond percolation}

\begin{aug}
\author[A]{\fnms{Geoffrey R.} \snm{Grimmett}\corref{}\thanksref{t1}\ead[label=e1]{g.r.grimmett@statslab.cam.ac.uk}\ead[label=u1,url]{http://www.statslab.cam.ac.uk/\textasciitilde grg/}}
\and
\author[A]{\fnms{Ioan} \snm{Manolescu}\thanksref{t2}\ead[label=e2]{i.manolescu@statslab.cam.ac.uk}}
\runauthor{G. R. Grimmett and I. Manolescu}
\affiliation{University of Cambridge}
\address[A]{Statistical Laboratory\\
Centre for Mathematical Sciences\\
University of Cambridge\\
Wilberforce Road\\
Cambridge CB3 0WB\\
United Kingdom\\
\printead{e1}\\
\hphantom{E-mail: }\printead*{e2}\\
\printead{u1}} %adresu isvedimo komanda gale!
\end{aug}

\thankstext{t1}{Supported in part by the EPSRC under Grant
EP/103372X/1.}

\thankstext{t2}{Supported by the EPSRC and Cambridge University.}

% HISTORY:
\received{\smonth{8} \syear{2011}}
\revised{\smonth{12} \syear{2011}}

% ABSTRACT
%
\begin{abstract}
All (in)homogeneous bond percolation models on the square, triangular,
and hexagonal lattices belong to the same universality class, in the
sense that they have identical critical exponents \textit{at} the
critical point (assuming the exponents exist). This is proved using the
star--triangle transformation and the box-crossing property. The
exponents in question are the
one-arm exponent~$\rho$, the $2j$-alternating-arms exponents
$\rho_{2j}$ for $j \ge1$, the volume exponent $\delta$, and the
connectivity exponent $\eta$. By earlier results of Kesten, this
implies universality also for the near-critical exponents $\beta$,
$\gamma$, $\nu$, $\Delta$ (assuming these exist) for any of these
models that satisfy a certain additional hypothesis, such as the
\textit{homogeneous} bond percolation models on these three lattices.
\end{abstract}

% KEYWORDS
%

\end{frontmatter}

%s1 #&#
\section{Introduction and results}
%s1.1 #&#
\subsection{Overview}
Two-dimensional percolation has enjoyed an extraordinary renaissance
since Smirnov's proof in 2001 of Cardy's formula (see~\cite{Smirnov}).
Remarkable progress has been made toward a full understanding of
site percolation on the triangular lattice, at and near its critical point.
Other critical two-dimensional models have, however, resisted solution.
The purpose of the current work
is to continue our study (beyond~\cite{GM1})
of the phase transition for inhomogeneous bond percolation on the
square, triangular and hexagonal lattices.
Our specific target is to show that such models belong to the same
universality class.
We prove that critical exponents \textit{at the critical point} are
constant within this
class of models (assuming that such exponents exist). We indicate a hypothesis
under which exponents \textit{near criticality} are constant also,
and note that
this is satisfied by the homogenous models.

We focus here on the one-arm exponent $\rho$, and the
$2j$-alternating-arms exponents $\rho_{2j}$ for $j \ge1$. By
transporting open primal paths and open dual paths,
we shall show that these exponents
are constant across (and beyond) the above class of bond percolation models.
More precisely, if any one of these exponents, $\pi$ say, exists for
one of
these models, then $\pi$ exists and is equal for every such model. No
progress is made here on the problem of existence
of exponents.

Kesten~\cite{Kes87a} showed that the exponents $\delta$ and $\eta$
are specified by knowledge
of $\rho$, under the hypothesis that $\rho$ exists. Therefore,
$\delta$ and $\eta$
are universal across this class of models.
Results related to those of~\cite{Kes87a} were obtained in~\cite{Kesten87} for
the ``near-critical'' exponents $\beta$, $\gamma$, $\nu$, $\De$. This
last work required
a condition of rotation-invariance not possessed by the
\textit{strictly} inhomogeneous models.
This is discussed further in Section~\ref{secres}.

It was shown
in~\cite{GM1} that critical inhomogeneous models on the above three
lattices have the box-crossing property;
this was proved by transportation of open box-crossings from the
homogeneous square-lattice model.
This box-crossing property, and the star--triangle transformation
employed to prove it, are the basic
ingredients that permit the
proof of universality presented here.

A different extension of the star--triangle method has been the subject of
work described in~\cite{BR,Z1,ZScull}. That work is, in a sense,
combinatorial in nature,
and it provides connections between percolation on a graph embedded in
$\RR^2$ and
on a type of dual graph obtained via
a generalized star--triangle transformation. In contrast, the work
reported here is closely
connected to the property of
isoradiality (see~\cite{BdTB,Ken02}), and is thus more geometric in
nature. It permits
the proof of relations \textit{between} a variety of two-dimensional graphs.
The connection to isoradiality will be the subject of a later paper~\cite{GM3}.

The paper is organized as follows. The relevant critical exponents are
summarized in Section~\ref{secexp},
and the main theorems stated in Section~\ref{secres}. Extensive
reference will be made to~\cite{GM1}, but the current work is fairly
self-contained.
Section~\ref{sectransformation} contains a short account of the
star--triangle transformation, for
more details of which the reader is referred to~\cite{GM1}. The proofs
are to be
found in Section~\ref{secexp2}.

%s1.2 #&#
\subsection{The models}\label{secmod}

Let $G = (V,E)$ be a countable connected planar graph,
embedded in $\RR^2$.
The bond percolation model on $G$ is defined as follows.
A~\textit{configuration} on
$G$ is an element $\omega=(\omega_e\dvtx e \in E)$ of the set $\Om=\{
0,1\}^E$.
An edge with endpoints $u$, $v$ is denoted $uv$.
The edge $e$ is called \textit{open} or $\omega$-\textit{open}
(resp., \textit{closed}) if $\omega_e=1$
(resp., $\omega_e = 0$).

For $\omega\in\Om$ and $A,B \subseteq V$, we say $A$ \textit{is
connected to} $B$ (in $\omega$),
written $A \leftrightarrow B$ (or $A\xleftrightarrow{G,\omega} B$),
if $G$ contains a path of open edges from some $a\in A$ to some $b \in B$.
An \textit{open cluster} of $\omega$ is a maximal set of
pairwise-connected vertices,
and the open cluster containing the vertex $v$ is denoted $C_v$.
We write $v \lra\oo$ if $v$ is the endpoint of an infinite open
self-avoiding path.\vadjust{\goodbreak}

The \textit{homogeneous} bond percolation model on $G$ is that
associated with the product measure
$\PP_p$ on $\Om$ with constant intensity $p\in[0,1]$. Let $0$ denote
a designated vertex
of $V$ called the \textit{origin}. The
\textit{percolation probability} and \textit{critical probability} are
given by
\begin{eqnarray*}
\theta(p) &=& \PP_p(0 \lra\infty),\\
\pc(G) &=& \sup\{p\dvtx\theta(p)=0\}.
\end{eqnarray*}

\begin{figure}

\includegraphics{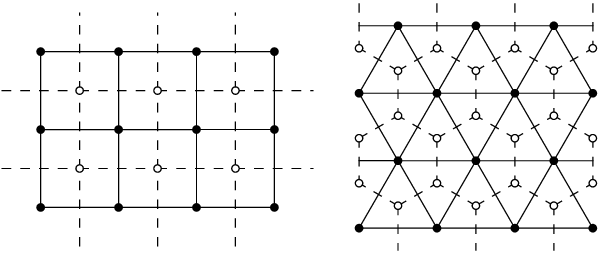}

\caption{The square lattice and its dual square lattice. The
triangular lattice and its dual hexagonal lattice.}
\label{figlatticeexamples}
\end{figure}

We consider the square, triangular,
and hexagonal (or honeycomb) lattices of Figure \ref
{figlatticeexamples}, denoted, respectively,
as $\ZZ^2$, $\TT$ and $\HH$.
It is standard that $\pc(\ZZ^2) = \frac12$, and $\pc(\TT)=1-\pc
(\HH)$ is the
root in the interval
$(0,1)$ of the cubic equation $3p-p^3-1=0$. See the references in
\cite{GrimmettPercolation,GM1} for
these and other known facts quoted in this paper.

%f1 #&#

We turn now to \textit{inhomogeneous} percolation on the above three lattices.
The edges of the square lattice are partitioned into two classes
(horizontal and vertical) of parallel edges,
while those of the triangular and hexagonal lattices may be split into
three such classes.
We allow the product measure on $\Om$ to
have different intensities on different edges, while requiring that any
two parallel edges have the same intensity.
Thus, inhomogeneous percolation on the square lattice has two parameters,
$p_0$ for horizontal edges and $p_1$ for vertical edges, and
we denote the corresponding measure $\PP_{\bp}^\squ$ where $\bp=(p_0,p_1)$.
On\vspace*{1pt} the triangular and hexagonal lattices, the measure is defined by a
triplet of parameters $\bp=(p_0, p_1, p_2)$,
and we denote these measures $\PP_{\bp}^\tri$ and $\PP_{\bp}^\hex
$, respectively.
Let $\sM$ denote the set of all such inhomogeneous bond percolation
models on the square, triangular, and hexagonal lattices,
with edge-parameters belonging to the half-open interval $[0,1)$.

These models have percolation probabilities and critical surfaces, and
the latter were
given explicitly in~\cite{GrimmettPercolation,GM1,Kestenbook}.
Let
\begin{eqnarray*}
\kappa_\squ(\bp) &=& p_0+p_1-1,\qquad \bp=(p_0,p_1), \\
\kappa_\tri(\bp) &=& p_0+p_1+p_2 - p_0 p_1 p_2 - 1,\qquad \bp
=(p_0,p_1,p_2),\\ %\label{criticaltriangular}\\
\kappa_\hex(\bp) &=& -\kappa_\tri(1-p_0,1-p_1,1-p_2),\qquad \bp
=(p_0,p_1,p_2).
\end{eqnarray*}
It is well known that the critical surface of the lattice $\ZZ^2$
(resp., $\TT$, $\HH$) is
given by $\kappa_\squ= 0$ [resp., $\kappa_\tri(\bp)=0$, $\kappa
_\hex(\bp)=0$].
Bond percolation on $\ZZ^2$ may be obtained from that on $\TT$ by
setting one parameter to zero.

The triplet $\bp=(p_0, p_1, p_2)\in[0,1)^3$ is called
\textit{self-dual} if
$\kappa_\tri(\bp)=0$.
We write $\alpha\pm\bp$ for the triplet $(\alpha\pm p_0,\alpha\pm
p_1,\alpha\pm p_2)$,
and also $\NN=\{1,2,\ldots\}$ for the natural numbers, and $\ZZ=\{
\ldots,-1,0,1,\ldots\}$
for the integers.

%s1.3 #&#
\subsection{Critical exponents}\label{secexp}

The percolation singularity is believed to be
of power-law type, and to be described by a number of
so-called ``critical exponents.'' These may be divided into two groups of
exponents: \textit{at criticality}
and \textit{near criticality}.

First, some notation: we write $f(t) \eqv g(t)$ as $t\to t_0 \in[0,\oo]$
if there exist strictly positive constants $A$, $B$ such that
\[
A g(t) \le f(t) \le Bg(t)
\]
in some neighbourhood of $t_0$ (or for all large $t$ in the case
$t_0=\oo$). We write $f(t) \logeqv g(t)$ if
$\log f(t) / \log g(t) \to1$. Two vectors $\bp_1=(p_1(e))$, $\bp
_2=(p_2(e))$ satisfy
$\bp_1<\bp_2$ if $p_1(e) \le p_2(e)$ for all $e$, and $\bp_1 \ne\bp_2$.

For simplicity, we restrict ourselves to the percolation models of the
last section.
Let $\LL=(V,E)$ be one of the square, triangular, and hexagonal lattices,
with origin denoted $0$.
Let $\bp= (p(e)\dvtx e\in E) \in[0,1)^E$ be invariant under translations
of $\LL$
as above, and let $\omega\in\Om$.
The lattice $\LL$ has a dual lattice $\LL^*=(V^*,E^*)$,
each edge of which is called \textit{open$^*$}
if it crosses a closed edge of~$\LL$. Open paths of $\LL$ are said to
have colour $1$,
and open$^*$ paths of $\LL^*$ colour~$0$.
We shall\vspace*{1pt} make use of duality as described in
\cite{GrimmettPercolation},
Section 11.2.

Let $\La_n$ be the set of all vertices within graph-theoretic
distance $n$ of the origin~$0$, with boundary $\pd\La_n = \La_n
\setminus\La_{n-1}$.
Let $\Ann(N,n) = \La_n \setminus\La_{N-1}$ be the \textit{annulus}
centred at $0$,
with \textit{interior radius} $N$ and \textit{exterior radius} $n$.
We call $\pd\La_n$ (resp., $\pd\La_N$) its \textit{exterior}
(resp., \textit{interior}) \textit{boundary}.
We shall soon consider embeddings of planar lattices in $\RR^2$,
and it will then be natural to use the $L^\oo$ metric rather than
graph-distance.
The choice of metric is in fact of no fundamental
important for what follows.
For $v \in V$, we write
\[
\rad(C_v) = \sup\{n\dvtx v \lra v+\pd\La_n\}.
\]

Let $\bpc$ be a vector lying on the critical surface. Thus, $\bpc$ is critical
in that
\[
\theta(\bp):= \PP_\bp(0\lra\oo)
\cases{ = 0, &\quad if $\bp< \bpc$,\cr
>0, &\quad if $\bp> \bpc$.}
\]
Let $k \in\NN$, and let
$\si=(\si_1,\si_2,\ldots,\si_k) \in\{0,1\}^k$; we call $\si$ a
\textit{colour sequence}. The sequence $\si$ is called \textit{monochromatic}
if either $\si=(0,0,\ldots,0)$ or $\si=(1,1,\ldots,1)$, and \textit
{bichromatic}
otherwise. If $k$ is even, $\si$ is called \textit{alternating} if
either $\si=(0,1,0,1,\ldots)$
or $\si=(1,0,1,0,\ldots)$.
For $0<N<n$,
the arm event $A_{\si}(N,n)$
is the event that the inner boundary of
$\Ann(N,n)$ is connected to its outer boundary by $k$ vertex-disjoint
paths
with colours $\sigma_1, \ldots, \sigma_k$, taken in anticlockwise order.
(Here and later, we require arms to be \textit{vertex}-disjoint rather
than \textit{edge}-disjoint.
This is an innocuous assumption since we work in this paper
with alternating colour sequences only.)

The choice of $N$ is in part immaterial to the study of the asymptotics of
$\PP_\bpc[A_{\si}(N,n)]$ as $n\to\oo$, and we shall assume
henceforth that $N=N(\si)$ is sufficiently
large that, for $n \ge N$, there exists a configuration with the
required $k$ coloured paths.
It is believed that there exist constants $\rho(\si)$ such that
\[
\PP_\bpc[\Arm_{\sigma}(N,n)] \logeqv n^{-\rho(\si)},
\]
and these are the \textit{arm-exponents} of the model. (Such asymptotics
are to be understood in the limit as $n\to\oo$.)

We concentrate here on the following exponents given in terms of $\PP
_\bpc$, with limits as $n\to\oo$:
\begin{longlist}[(a)]
\item[(a)]
volume exponent: $\PP_\bpc(|C_0| = n) \logeqv n^{-1-1/\delta}$,
\item[(b)]
connectivity exponent: $\PP_\bpc(0 \lra x) \logeqv|x|^{-\eta}$,
\item[(c)]
one-arm exponent: $\PP_\bpc(\rad(C_0) = n) \logeqv n^{-1-1/\rho}$,
\item[(d)]
$2j$-alternating-arms exponents: $\PP_\bpc[\Arm_{\sigma}(N,n)]
\logeqv n^{-\rho_{2j}}$,
for each alternating colour sequence $\si$ of length $2j$.
\end{longlist}
It is believed that the above asymptotic relations hold for
suitable exponent-values, and indeed with $\logeqv$
replaced by the stronger relation $\eqv$. Essentially the only two-dimensional
percolation process for which these limits are proved
(and, furthermore, the exponents calculated explicitly)
is site percolation on the triangular lattice (see
\cite{Smirnov,Smirnov-Werner}).

The arm events are defined above in terms of open primal and open$^*$ dual
paths. When considering site percolation, one considers instead
open paths in the primal and matching lattices. This is especially simple
for the triangular lattice since $\TT$ is self-matching.
It is known for site percolation on the triangular lattice,~\cite{ADA},
that for given $k \in\NN$, the exponent
for $\rho(\si)$ is constant for any bichromatic
colour sequence $\si$ of given length $k$. This is believed to
hold for other two-dimensional models also, but no proof is known.
In particular, it is believed for any model in $\sM$ that
\[
\PP_\bpc[\Arm_{\sigma}(N,n)] \logeqv n^{-\rho_{2j}}
\]
for any bichromatic colour sequence $\si$ of length $2j$, and any $j
\ge1$.

We turn now to the near-critical exponents which, for definiteness, we
define as follows.
Let $\bp=(p(e)\dvtx e \in E) \in[0,1)^E$ and $\eps\in\RR$,
and write $\PP_{\bp+\eps}$ for the product measure on $\Om$ in
which edge $e$ is open
with probability
\[
(\bp+\eps)_e:= \max\bigl\{0,\min\{p(e)+\eps,1\}\bigr\}.
\]
By subcritical exponential-decay
(see~\cite{GrimmettPercolation}, Section 5.2), for $\eps>0$, there exists
$\xi= \xi(\bpc-\eps)\in[0,\oo)$ such that
\[
-\frac1n \log\PP_{\bpc-\eps}(0\lra\pd\La_n) \to1/\xi\qquad
\mbox{as } n \to\oo.
\]
The function $\xi$ is termed the \textit{correlation length}.

Here are the further exponents considered here:
\begin{longlist}[(a)]
\item[(a)]
percolation probability: $\theta(\bpc+\eps) \logeqv\eps^\beta$ as
$\eps\downarrow0$,
\item[(b)]
correlation length: $\xi(\bpc-\eps) \logeqv\eps^{-\nu}$ as $\eps
\downarrow0$,
\item[(c)]
mean cluster-size: $\PP_{\bpc+\eps}(|C_0|; |C_0|<\oo) \logeqv|\eps
|^{-\gamma}$ as $\eps\to0$,
\item[(d)]
gap exponent: for $k \ge1$, as $\eps\to0$,
\[
\frac{\PP_{\bpc+\eps}( |C_0|^{k+1}; |C_0|<\oo)}{\PP_{\bpc+\eps
}( |C_0|^{k}; |C_0|<\oo)}
\logeqv|\eps|^{-\De}.
\]
\end{longlist}
We have written $\PP(X)$ for the mean of $X$ under the probability
measure $\PP$, and
$\PP(X;A)=\PP(X1_A)$ where $1_A$ is the indicator function of the
event $A$.

As above, the near-critical exponents
are known to exist essentially only for site percolation on the
triangular lattice.
See~\cite{GrimmettPercolation}, Chapter 9, for a general account
of critical exponents and scaling theory.

%s1.4 #&#
\subsection{Principal results}\label{secres}

A critical exponent $\pi$ is said to \textit{exist} for a model $M\in
\sM$ if the
appropriate asymptotic relation (above) holds, and $\pi$ is called
\textit{$\sM$-invariant} if it exists for all $M \in\sM$ and its value
is independent of the choice of such $M$.
%
%th1.1 #&#
\begin{theorem}\label{thmeq}
For every $\pi\in\{\rho\}\cup\{\rho_{2j}\dvtx j \ge1\}$, if $\pi$
exists for some model $M \in\sM$, then
it is $\sM$-invariant.
\end{theorem}

By the box-crossing property of~\cite{GM1}, Theorem 1.3, we may apply
the theorem of
Kesten~\cite{Kes87a}
to deduce the following. If either $\rho$ or $\eta$ exists for some
$M \in\sM$, then:
\begin{longlist}[(a)]
\item[(a)]
both $\rho$ and $\eta$ exist for $M$,
\item[(b)] $\delta$ exists for $M$,
\item[(c)]
the scaling relations $\eta\rho=2$ and $2\rho= \delta+1$ are valid.
\end{longlist}
Taken in conjunction with Theorem~\ref{thmeq}, this implies in particular
that $\delta$ and $\eta$ are $\sM$-invariant
whenever either $\rho$ or $\eta$ exist for some $M\in\sM$.

\begin{figure}

\includegraphics{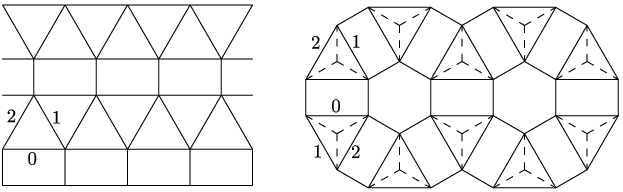}

\caption{Isoradial embeddings of the Archimedean lattices $(3^3,4^2)$
and $(3,4,6,4)$.
The second may be
transformed into the hexagonal lattice by one sequence of
star--triangle transformations as marked.
An edge parallel to one labelled $i$ has edge-parameter $p_i$,
and the two parameters on any square have sum $1$.}
\label{figarch1}
\end{figure}

We note in passing that Theorem~\ref{thmeq} may be extended to
certain other graphs derived from
the three main lattices of this paper by sequences of star--triangle
transformations, as well
as to their dual
graphs. This includes a number of
tessellations (see~\cite{GS}) and, in particular, two further
Archimedean lattices,
namely those denoted $(3^3,4^2)$ and $(3,4,6,4)$ and illustrated
in Figure~\ref{figarch1}. The measures on these two lattices are as
follows. Let $\bp=(p_0,p_1,p_2)\in[0,1)^3$
be self-dual. Edge $e$ is open with probability $p(e)$ where:
\begin{longlist}[(a)]
\item[(a)] $p(e) = p_0$ if $e$ is horizontal,
\item[(b)] $p(e) = p_1$ if $e$ is parallel to the right edge of an upward
pointing triangle,
\item[(c)] $p(e)=p_2$ if $e$ is parallel to the left edge of an upward
pointing triangle,
% \item$p(e) = 1-p_2$ if $e$ is the right edge of an upward pointing
%star,
% \item$p(e)=1-p_1$ if $e$ is the left edge of an upward pointing
%star.
%
\item[(d)] the two parameters of any rectangle have sum $1$.
\end{longlist}
Theorem~\ref{thmeq} holds with $\sM$ augmented by all such bond
models on
these two lattices. The proofs are essentially the same.
The methods used here do not appear to extend to homogeneous percolation
on these two lattices.
Drawings of the eleven Archimedean lattices and their duals may be found
in~\cite{PW}.
We note that the remaining six Archimedean lattices
may not be embedded isoradially in the plane (see~\cite{GM3,KenS}).

%f2 #&#

There is a simple reason for the fact that
Theorem~\ref{thmeq} concerns the alternating-arm exponents rather
than all arm exponents.
We shall see in Section
\ref{sectransformation} that the star--triangle transformation
conserves open primal and
open$^*$ dual paths,
but that, in certain circumstances,
it allows distinct paths
of the same colour to coalesce.

The box-crossing property of~\cite{GM1}
implies an affine isotropy of these models at criticality, yielding in
particular that certain directional exponents
are independent of the choice of direction. For example, let $\theta
\in[0,\pi)$, and
consider the probability of an open path from the origin
to a line with gradient $\tan\theta$ and distance $\pm n$ from the origin.
The associated exponent equals the undirected exponent $\rho$.
A similar
statement holds for arm-directions in the alternating-arm exponents.
These facts follow by the box-crossing property (in conjunction with
the separation theorem of Section~\ref{secsep}), in particular
by its consequence that any annulus comprising rectangles of given aspect-ratio
contains an open cycle with probability bounded away from $0$.

Kesten has shown in~\cite{Kesten87} (see also~\cite{Nolin})
that the above near-critical exponents may be
given explicitly in terms of exponents at criticality, for
two-dimensional models
satisfying certain hypotheses. Homogeneous percolation on our three lattices
satisfy these hypotheses, but it is not known whether the strictly inhomogeneous
models have sufficient regularity for the conclusions to hold for them.
The basic problem is that,
while the box-crossing property of~\cite{GM1} implies an isotropy for
these models at
criticality,
the corresponding isotropy away from criticality is unknown.
For this reason, we restrict the statement of the next theorem to homogeneous
models.
%
%th1.2 #&#
\begin{theorem}\label{thmeq2}
Assume that $\rho$ and $\rho_4$ exist for some $M \in\sM$.
Then $\beta$, $\gamma$, $\nu$, and $\De$ exist
for homogeneous percolation on the square, triangular and hexagonal
lattices, and they are invariant
across these three models. Furthermore, they satisfy the scaling relations
\[
\rho\beta=\nu,\qquad \rho\gamma=\nu(\delta-1),\qquad \rho\De=\nu
\delta.
\]
\end{theorem}

The proof is an adaptation of the arguments and conclusions of
\cite{Kesten87,Nolin}, and is
omitted here.

Other authors have observed hints of universality, and we mention for
example~\cite{SedW},
where it is proved that certain dual pairs of lattices have equal
exponents (whenever these exist).

%f3 #&#
\begin{figure}[b]

\includegraphics{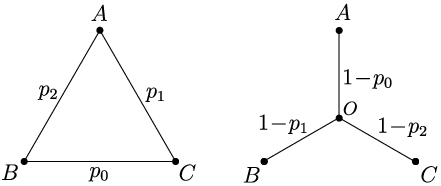}

\caption{The star--triangle transformation when $\kappa_\tri(\bp)=0$.}
\label{figstartriangletransformation}
\end{figure}

%s2 #&#
\section{Star--triangle transformation}\label{sectransformation}

Consider the triangle $G=(V,E)$
and the star $G'=(V',E')$, as drawn in
Figure~\ref{figstartriangletransformation}. Let $\bp
=(p_0,p_1,p_2)\in[0,1)^3$.
Write $\Om=\{0,1\}^E$ with associated product probability measure $\PP
_\bp^\tri$,
and $\Om'=\{0,1\}^{E'}$ with associated measure $\PP_{1-\bp}^\hex$,
as illustrated
in the figure.
Let $\omega\in\Om$ and
$\omega'\in\Om'$.
For each graph, we may consider open connections between its vertices,
and we abuse notation by writing, for example, $x \xleftrightarrow
{G,\omega} y$ for the \textit{indicator
function} of the event that $x$ and $y$ are connected\vadjust{\goodbreak} by an open path
of~$\omega$.
Thus, connections in $G$ are described by the family
$( x \xleftrightarrow{G, \omega} y\dvtx x,y \in V)$ of random variables,
and similarly for $G'$.

It may be shown that the two families
\[
\bigl( x \xleftrightarrow{G, \omega} y\dvtx x,y = A,B,C\bigr),\qquad
\bigl( x \xleftrightarrow{G', \omega'} y\dvtx x,y = A,B,C\bigr),
\]
of random variables have the same joint law whenever $\kappa_\tri(\bp)=0$.
That is to say, if $\bp$ is self-dual, the existence (or not) of open
connections
is preserved (in law) under the star--triangle transformation. See
\cite{GrimmettPercolation}, Section 11.9.

\begin{figure}

\includegraphics{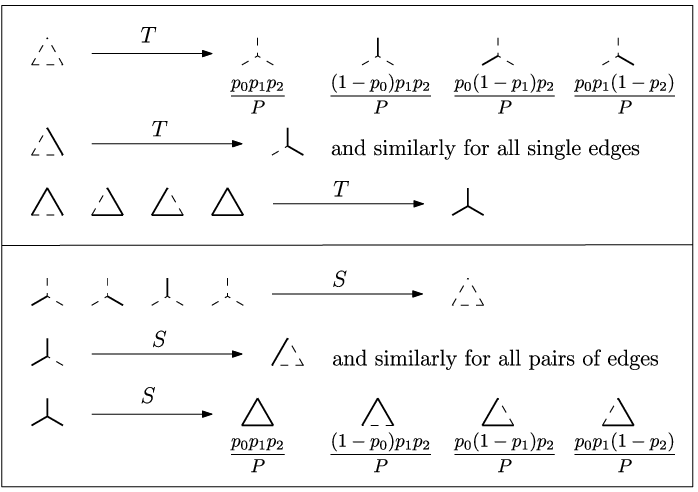}

\caption{The ``kernels'' $T$ and $S$ and their
transition probabilities. The constant $P$
is given by $P:=(1-p_0)(1-p_1)(1-p_2)$.}
\label{figsimpletransformationcoupling}
\end{figure}

The two measures $\PP_\bp^\tri$ and $\PP_{1-\bp}^\hex$ may be
coupled in a natural
way. Let $\bp\in[0,1)^3$ be self-dual, and let $\Om$ (resp., $\Om'$)
have associated measure $\PP_\bp^\tri$ (resp., $\PP_{1-\bp}^\hex$)
as above. The random mappings $T\dvtx\Om\to\Om'$ and
$S\dvtx\Om'\to\Om$ of Figure~\ref{figsimpletransformationcoupling}
are such that: $T(\omega)$ has law $\PP_{1-\bp}^\hex$, and
$S(\omega')$ has law $\PP_\bp^\tri$. Under this coupling,
the presence or absence of connections between the corners $A$, $B$,
$C$ is preserved.

%f4 #&#

\begin{figure}

\includegraphics{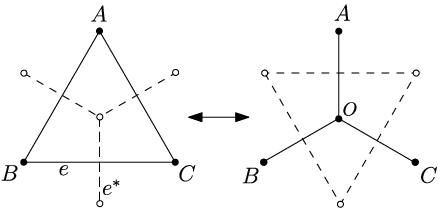}

\caption{The star--triangle transformation acts simultaneously on
primal and dual graph elements.}
\label{figdual}
\end{figure}

The maps $S$ and $T$ act on configurations on stars and triangles. They
act simultaneously on the
duals of these graph elements, illustrated in Figure~\ref{figdual}.
Let $\omega\in\Om$, and define $\omega^*(e^*)=1-\omega(e)$
for each primal/dual pair $e/e^*$ of the left-hand side of the
figure. The
action of $T$ on $\Om$ induces an action on
the dual space $\Om^*$, and it is easily checked that this action
preserves $\omega^*$-connections. The map
$S$ behaves similarly. This property of the star--triangle
transformation has been generalized
and studied in~\cite{BR} and the references therein.

%f5 #&#

%f6 #&#
\begin{figure}[b]

\includegraphics{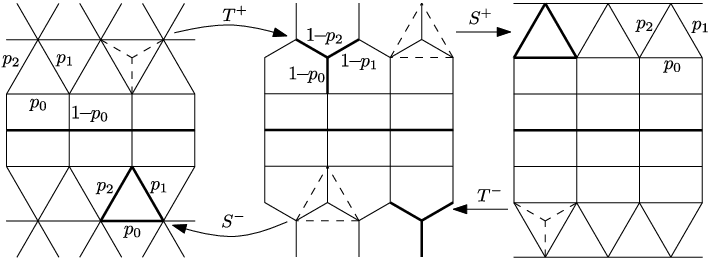}

\caption{The transformation
$S^+ \comp T^+$ (resp., $S^- \comp T^-$)
transforms $\Lattice^1$ into $\Lattice^2$
(resp., $\Lattice^2$ into $\Lattice^1$).
They map the dashed graphs to
the bold graphs.}
\label{figlatticestriptransformation}
\end{figure}

So-called mixed lattices were introduced in~\cite{GM1}. These are
hybrid embeddings of
the square lattice with either the triangular or hexagonal lattice, the
two parts being
separated by a horizontal interface. By means of appropriate
star--triangle transformations,
the interface may be moved up or down, and this operation permits the
transportation of
open box-crossings between the square lattice and the other lattice.
Whereas this was
suited for proving the box-crossing property, a slightly altered hybrid is
useful for studying arm exponents.

Let $m \ge0$, and consider the mixed lattice $\Lattice^m=(V^m,E^m)$
drawn on the left of Figure~\ref{figlatticestriptransformation},
formed of a horizontal strip of the square lattice
centred on the $x$ axis of height $2m$, with the triangular lattice
above and beneath it.
The embedding of each lattice is otherwise as in~\cite{GM1}:
the triangular lattice comprises equilateral triangles with side length
$\sqrt{3}$,\vspace*{1pt}
and the square lattice comprises rectangles
with horizontal dimension
$\sqrt{3}$ and vertical dimension~$1$.
We require also that the origin of $\RR^2$ be a vertex of the mixed lattice.

Let\vspace*{1pt} $\bp\in[0,1)^3$, and let $\PP_\bp^m$ be the product measure on
$\Om^m = \{0,1\}^{E^m}$
for which edge $e$ is open with probability $p(e)$ given by:
\begin{longlist}[(a)]
\item[(a)] $p(e) = p_0$ if $e$ is horizontal,
\item[(b)] $p(e) = 1- p_0$ if $e$ is vertical,
\item[(c)] $p(e) = p_1$ if $e$ is the right edge of an upward pointing
triangle,
\item[(d)] $p(e)=p_2$ if $e$ is the left edge of an upward pointing triangle.
% \item$p(e) = 1-p_2$ if $e$ is the right edge of an upward pointing
%star,
% \item$p(e)=1-p_1$ if $e$ is the left edge of an upward pointing
%star.
\end{longlist}

Suppose further that $\bp$ is self-dual, in that $\kappa_\tri(\bp
)=0$, and let $\omega^m\in\Om^m$.
We denote by $\Tu$ (resp., $\Td$) the transformation $T$ of Figure
\ref{figsimpletransformationcoupling}
applied to an upward (resp., downward) pointing triangle.
Write $T^+$ for the transformation of $\omega$ obtained by applying
$\Td$ to every downward
pointing triangle in the upper half plane,
and $\Tu$ similarly in the lower half plane; sequential
applications of star--triangle transformations are required to be
independent of one another.

Similarly, we denote by $\Su$ (resp., $\Sd$) the transformation $S$
of Figure~\ref{figsimpletransformationcoupling}
applied to an upward (resp., downward) pointing star.
Write
$S^+$ for the transformation of $(T^+ \Lattice^m, T^+(\omega^m))$
obtained by applying $\Su$ to all upward pointing stars in the upper
half-plane
and similarly $\Sd$ in the lower half-plane.
It may be checked that $\omega^{m+1} = S^+\circ T^+(\omega^m)$ lies
in $\Om^{m+1}$
and has law $\PP_\bp^{m+1}$. That is, viewed as a transformation
acting on measures,
we have $(S^+\circ T^+) \PP_\bp^m = \PP_\bp^{m+1}$.

The transformations $T^-$ and $S^-$ are defined similarly, and illustrated
in Figure~\ref{figlatticestriptransformation}. As in that figure,
for $m \ge0$,
\begin{eqnarray*}
(S^+ \comp T^+) \Lattice^m &=& \Lattice^{m+1},\qquad
(S^+ \comp T^+) \PP_\bp^m = \PP_\bp^{m+1},\\
(S^- \comp T^-) \Lattice^{m+1} &=& \Lattice^{m},\qquad
(S^- \comp T^-) \PP_\bp^{m+1} = \PP_\bp^m.
\end{eqnarray*}

We turn to the operation of these two transformations on open paths,
and will concentrate on
$S^+ \circ T^+$; similar statements are valid for $S^-\circ T^-$. Let
$\omega^m\in\Om^m$, and
let $\pi$ be an $\omega^m$-open path of $\Lat^m$. It is not
difficult to see (and is explained fully in
\cite{GM1}) that the image of $\pi$ under $S^+\circ T^+$ contains
some $\omega^{m+1}$-open path $\pi'$.
Furthermore, $\pi'$ lies within the $1$-neighbourhood of $\pi$ viewed
as a subset of $\RR^2$, and has endpoints within unit Euclidean
distance of
those of $\pi$. Any vertex of $\pi$ in the square part of $\Lat^m$
is unchanged
by the transformation. The corresponding statements hold also for
open$^*$ paths of the dual of $\Lat^m$.
These facts will be useful in observing the effect of $S^+ \circ T^+$
on the arm events.

Arm exponents are defined in Section~\ref{secexp} in terms of boxes
that are adapted
to the lattice viewed as a graph. It will be convenient to work also
with boxes of $\RR^2$.
Let $\Lat=(V,E)$ be a mixed lattice duly embedded in~$\RR^2$, and
write $V_0$
for the subset of $V$ lying on the $x$-axis.
Let $\omega\in\Om=\{0,1\}^E$. For $R \subseteq\RR^2$ and $A,B
\subseteq R \cap V_0$,
we write $A \xlra{R,\omega} B$ (with negation written $A \nxlra
{R,\omega} B$) % (or $A \lra B$ in $R/\om$)
if there\vadjust{\goodbreak} exists an $\omega$-open path joining some $a\in A$ and some
$b \in B$
using only edges that intersect $R$. Let $D$ be the unit (Euclidean)
disk of $\RR^2$ and write
$R + D$ for the direct sum $\{r+d\dvtx r\in R, d\in D\}$.
%
%pr2.1 #&#
\begin{prop}\label{proppaths}
Let $m \ge0$, $\omega\in\Om^m$, $R \subseteq\RR^2$, and $u,v \in
R\cap V_0$.
For $\tau\in\{S^+ \comp T^+, S^- \comp T^-\}$:
\begin{longlist}[(a)]
\item[(a)] if $u \xlra{R,\omega} v$ then
$u\xlra{R+D, \tau(\omega)} v$,
\item[(b)] if $u \nxlra{R+D,\omega} v$
then $u \nxlra{R,\tau(\omega)} v$.
\end{longlist}
\end{prop}
\begin{pf}
(a) Let $\tau=S^+ \comp T^+$; the case $\tau=S^- \comp T^-$
is similar (we assume $m \ge1$ where necessary).
If $u \xlra{R,\omega} v$, there exists an $\omega$-open path $\pi$
of $\Lat$ from $u$ to $v$ using
edges that intersect $R$.
Since $u$, $v$ are not moved by $\tau$, the image $\tau(\pi)$
contains a $\tau(\omega)$-open path
of $\tau\Lat$ from $u$ to $v$. It is elementary that $\tau$
transports paths through a distance not exceeding $1$ (see~\cite{GM1}, Proposition~2.4).
Therefore, every edge of $\tau(\pi)$ intersects $R+D$.

(b) Suppose\vspace*{1pt} $u \xlra{R,\tau(\omega)} v$. By considering the
star--triangle transformations
that constitute
the mapping $\tau$ [as in part (a)], we have that \mbox{$u \xlra{R+D,\omega} v$}.
\end{pf}

%s3 #&#
\section{Universality of arm exponents}\label{secexp2}

This section contains the proof of Theorem~\ref{thmeq}.
The reader is reminded that we work with translation-invariant measures
associated
with the square, triangular and hexagonal lattices.

%s3.1 #&#
\subsection{The arm exponents}\label{secexist}

Let $k \in\NN$ and $\si\in\{0,1\}^k$.
The arm event $\Arm_{\si}(N,n)$ is empty if $N$ is too small
to support the existence of the required
$k$ disjoint paths to the exterior boundary of the annulus $\Ann(N,n)$.
As explained in~\cite{Nolin}, for example, for each $\si$, there
exists $N=N_0(\si)$ such that the
arm exponent (assuming existence) is independent of the choice of $N
\ge N_0(\si)$. We assume henceforth
that $N$ is chosen sufficiently large for this to be the case.

It is a significant open problem of probability theory to prove the existence
and invariance of arm exponents for general lattices. This amounts to
the following in the present situation.
%
%co3.1 #&#
\begin{conj}\label{defarmexp}
Let $\bp\in[0,1)^3$ be self-dual.
For $k \in\NN$ and a colour sequence $\sigma\in\{0,1\}^k$,
there exists $\rho= \rho(\sigma,\bp) >0$ such that
\[
\PP_\bp^\tri[\Arm_{\sigma}(N,n)] \logeqv n^{-\rho}.
\]
Furthermore, $\rho(\si,\bp)$ is constant for all self-dual $\bp$.
\end{conj}

This is phrased for the triangular lattice, but it embraces also the
square and
hexagonal lattices, the first by setting a component of $\bp$ to $0$,
and the second
by a single application of the star--triangle transformation.\vadjust{\goodbreak} (See also~\cite{SedW}.)
We make no contribution toward a proof of the
claim of existence in this conjecture.
Theorem~\ref{thmeq} amounts to the proof of the claim of constantness,
for an alternating sequence $\si$ of length $k\in\{1,2,4,\ldots\}$.

Hereafter, we consider only the one-arm event with $\si=\{1\}$,
and the $2j$-alternating-arms events with $\si=(1,0,1,0,\ldots)$,
with associated exponents denoted, respectively, as $\rho_1$ and $\rho_{2j}$.
Thus, $\rho_1=1/\rho$ with $\rho$ as in Section~\ref{secexp}.

%s3.2 #&#
\subsection{The arm events}\label{secmainp}

Let $\Lat$ be one of the square, triangular, and hexagonal lattices,
or a hybrid
thereof as in Section~\ref{sectransformation}. We embed $\Lat$
in $\RR^2$ in the manner described in that section.
Let $x_i=(i\sqrt3,0)$, $i \ge0$, denote the
vertices common to these lattices to the right of the origin,
and $y_i=((i+\frac12)\sqrt3,\frac12)$, $i \ge0$, the vertices of
the dual lattice $\Lat^*$
corresponding to the faces of $\Lat$
lying immediately above the edge $x_i x_{i+1}$.
For $r\in(0,\oo)$, let $B_r=[-r,r]^2
\subseteq\RR^2$, with boundary $\pd B_r$. We recall that $C_x$
(resp., $C_y^*$)
denotes the open cluster of $\Lat$ containing $x$
(resp., the open$^*$ cluster of $\Lat^*$ containing $y$). For $n \ge
1$ and any connected subgraph
$C$ of either $\Lat$ or $\Lat^*$, we write $C \cap\pd B_r \ne\es$
if $C$ contains vertices in both
$B_r$ and $\RR^2 \setminus(-r,r)^2$.
Note that we may have $C \cap\pd B_r \ne\es$ even when no
vertices of $C$ lie in $\pd B_r$.

For $j,n\in\NN$ with $j \ge2$, let
\begin{eqnarray*}
\Arm_{1}(n) &=& \{ C_{x_0} \cap\pd B_n \ne\es\}, \\
\Arm_{2}(n) &=& \{ C_{x_0}\cap\pd B_n \ne\es, C^*_{y_0} \cap\pd
B_n \ne\es\}, \\
\Arm_{2j}(n) &=& \bigcap_{0\le i < j}
\bigl\{C_{x_{i}}\cap\pd B_n \ne\es\mbox{, and } x_i \nxlra
{B_n,\omega} \{x_0, x_1,\ldots,x_{i-1}\}\bigr\}.
\end{eqnarray*}
We write $\Arm_k^\Lat(n)$ when the role of $\Lat$ is to be stressed.
Note the condition of disconnection in the definition of $A_{2j}(n)$:
it is required that the $x_i$ are not connected by open paths
of edges all of which intersect $B_n$.

An alternative proof of the second inequality of the next lemma may be
obtained with the help of
the forthcoming separation theorem, Theorem~\ref{separation}, as in the
final part of the proof of Proposition~\ref{armextension}.
The latter route is more general since it assumes less about the
underlying lattice,
but it is also more complex since it relies on a version of
the separation theorem of~\cite{Kesten87} whose
somewhat complicated proof is omitted from the current work.
%
%le3.2 #&#
\begin{lemma}\label{lemhook}
Let $\bp\in[0,1)^E$ be self-dual. Let $k\in\{1,2,4,6,\ldots\}$, and
let $\si$ be an alternating colour sequence
of length $k$ (when $k=1$ we set $\si=\{1\}$).
There exists $N_0 =N_0(k)\in\NN$ and $c=c(\bp,N,k)>0$ such that
%
%e3.1 #&#
\begin{equation}\label{G300}
\PP_\bp^\tri\bigl[\Arm_k\bigl(n\sqrt3\bigr)\bigr] \le\PP_\bp^\tri[\Arm_\si(N,n)]
\le c\PP_\bp^\tri[\Arm_k(n)]
\end{equation}
for $n \ge N\ge N_0$.
\end{lemma}
\begin{pf}%{Outline proof of Lemma~\ref{lemhook}}
First, here is a note concerning the event $\Arm_{2j}(n)$ with $j \ge
2$. If $\omega\in\Arm_{2j}(n)$,
vertices\vadjust{\goodbreak} $x_i$, $0\le i < j$, are connected to $\pd B_n$ by open paths.
We claim that
$j$ such open paths may be found that are vertex-disjoint and interspersed
by $j$ open$^*$ paths joining the $y_i$ to $\pd B_n$. This will imply
the existence of
$2j$ arms of
alternating types joining $\{x_0,y_0,x_1,y_1,\ldots,x_{j-1}\}$ to
$\pd B_n$, such that the open primal paths are vertex-disjoint, and the
open$^*$ dual paths are vertex-disjoint except at the $y_i$.
The claim may be seen as follows (see also Figure~\ref{figconn}).
The dual edge $e$ with endpoints
$y_0$, $y_0-(0,1)$ is necessarily open$^*$.
By exploring the boundary of $C_{x_0}$ at $e$, one may find two
open$^*$ paths
denoted $\pi_0$, $\pi_0'$, joining $y_0$ to $\pd B_n$, and
vertex-disjoint except at $y_0$. Let $1\le r \le j-2$.
Since $x_r,x_{r+1} \xlra{B_n,\omega} \pd B_n$
and $x_r \nxlra{B_n,\omega} x_{r+1}$, we may similarly explore the
boundary of $C_{x_r}$ to
find an open$^*$ path $\pi_r$ of $B_n$ that joins $y_r$ to $\pd B_n$,
and is vertex-disjoint
from either $\pi_0$ or $\pi_0'$, and in addition from $\pi_s$, $s < r$.
The dual paths $\pi_0',\pi_0,\pi_1,\ldots,\pi_{j-2}$ are the required
open$^*$ arms.

The set $\La_n$ induces a subgraph of $\TT$ whose boundary is denoted
$\pd\La_n$.
We denote the inside of
$\pd\La_n$ (i.e., the closure of the bounded component of
$\RR^2 \setminus\pd\La_n$) by $\La_n$ also.
It is easily seen that $\La_n \subseteq B_{n\sqrt3}$,
and the first inequality of (\ref{G300}) follows immediately.

For the second inequality,
we shall use the fact that $B_n \subseteq\La_n$, together with
a suitable construction of open and open$^*$ paths within $\La_N$.
Let $k =2j \in\{2,4,6,\ldots\}$ and suppose $\Arm_\si(N,n)$ occurs.
On an anticlockwise traverse of
$\pd\La_N$, we find points $a_0,b_0,a_1,b_1,\ldots,a_{j-1},b_{j-1}$
such that the $a_i$ (resp., $b_i$) are endpoints of open (resp., open$^*$)
paths crossing the annulus $\Ann(N,n)$. Note that the $b_i$ are
not vertices of $\Lat^*$, but instead lie in open$^*$ edges.
Write $\mathbf a=(a_0,a_1,\ldots,a_{j-1})$,
$\mathbf b = (b_0,b_1,\ldots,b_{j-1})$.

%f7 #&#
\begin{figure}

\includegraphics{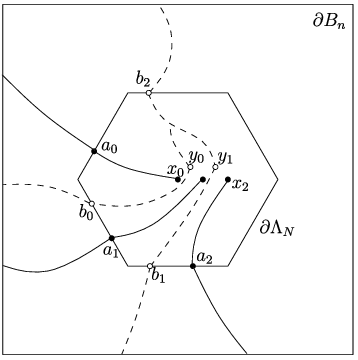}

\caption{The vertices $x_r$ (resp., $y_r$) may be connected by open
(resp., open$^*$)
paths to the $a_r$ (resp., $b_r$) in such a way that the event $\Arm
_6(n)$ results.}
\label{figconn}
\end{figure}

As illustrated in Figure~\ref{figconn}, for sufficiently large $N$
and all vectors $\mathbf a$, $\mathbf b$ of length~$j$,
there exists a configuration $\omega_{\mathbf a, \mathbf b}$
of primal edges of $\La_N$ such that
$x_r \xlra{\La_N} a_r$ for $0\le r \le j-1$, and $y_r \,{\xlra{\La_N}
\!\!\!{}^*} \,\,b_r$
and
$y_r \,{\xlra{\La_N}\!\!\!{}^*}\,\, b_{j+1}$ for $0\le r \le j-2$.
That is, conditional on $\Arm_\si(N,n)$, if $\omega_{\mathbf a,
\mathbf b}$ occurs then so
does $\Arm_k(n)$.
Assume for the moment that $p_i>0$ for all $i$.
The configurations $\omega_{\mathbf a, \mathbf b}$ may be chosen in
such a way that
\[
c'=c'(\bp, N,k):= \min_{\mathbf a,\mathbf b,\omega}
\PP_\bp^\tri(\omega_{\mathbf a, \mathbf b}\mid\Arm_\si(N,n))
\]
satisfies $c'>0$. Note that $c'$ does not depend on $n$.
The details of the construction of the $\omega_{\mathbf a, \mathbf b}$ are
slightly complicated but follow standard lines and are omitted
(similar arguments are used in~\cite{GrimmettPercolation}, Section
8.2, and
\cite{WernerSMF}, Chapter~2).
It follows as required that
\[
c' \PP_\bp^\tri[A_\si(N,n)] \le\PP_\bp^\tri[A_k(n)].
\]

Whereas a naive construction of the
$\omega_{\mathbf a, \mathbf b}$ succeeds when $p_i>0$
for all $i$,
a minor variant of the argument is needed if $p_i=0$ for some $i$.
The details are elementary and are omitted.

The case $k=1$ is similar but simpler.
\end{pf}

%s3.3 #&#
\subsection{\texorpdfstring{Proof of Theorem \protect\ref{thmeq}}{Proof of Theorem 1.1}}\label{secpfthm}

Let $\bp\in[0,1)^3$ be self-dual with $p_0>0$, and consider
the two measures $\PP_{(p_0, 1-p_0)}^\square$ (resp., $\PP_{\bp
}^\triangle$) on the
square (resp., triangular) lattice.
The proof of the universality of the box-crossing property was based on
a technique
that transforms one of these lattices into the other while preserving primal
and dual connections.
The same technique will be used here to prove the following result, the
proof of
which is deferred to Section~\ref{secpfarm}.
%
%pr3.3 #&#
\begin{prop}\label{armineq}
For any $k \in\{1, 2, 4, 6, \ldots\}$ and any self-dual triplet $\bp
\in[0,1)^3$ with $p_0>0$,
there exist $c_0, c_1, n_0 > 0$ such that, for all $n \geq n_0$,
\[
c_0 \PP_{\bp}^\triangle[\Arm_k(n)]
\leq\PP_{(p_0, 1-p_0)}^\square[\Arm_k(n)]
\leq c_1 \PP_{\bp}^\triangle[\Arm_k(n)].
\]
\end{prop}
%
% We will see in the proof of the proposition that the constants
%actually are
% uniform in all $\bp$ such that $\max_i \{p_i \} < 1-\varepsilon$.
%
\begin{pf*}{Proof of Theorem~\ref{thmeq}}
Suppose there exist $k \in\{1,2,4,6,\ldots\}$, a self-dual $\bp\in[0,1)^3$,
and $\alpha>0$, such that
%
%e3.2 #&#
\begin{equation}\label{G21}
\PP_\bp^\tri[\Arm_\si(N,n)] \logeqv n^{-\alpha}
\end{equation}
with $\si$ the alternating colour sequence of length $k$ (when $k=1$,
we take
$\si=\{1\}$).
By Lemma~\ref{lemhook}, (\ref{G21}) is equivalent to
%
%e3.3 #&#
\begin{equation}\label{G2}
\PP_\bp^\tri[\Arm_k(n)] \logeqv n^{-\alpha}.
\end{equation}
We say that ``$\PP$ satisfies (\ref{G2})''
if (\ref{G2}) holds with $\PP_\bp^\tri$ replaced
by $\PP$.

By self-duality, there exists $i$
such that $p_i>0$, and we assume without loss
of generality that $p_0 >0$. By Proposition~\ref{armineq},
$\PP_{(p_0,1-p_0)}^\squ$ satisfies (\ref{G2}).
Similarly, $\PP_{\bp'}^\tri$ satisfies (\ref{G2}) for any self-dual
$\bp'\in[0,1)^3$
of the form $\bp'=(p_0,p_1',p_2')$. The claim is proved after further
applications of
the proposition.
\end{pf*}

%s3.4 #&#
\subsection{\texorpdfstring{Proof of Proposition \protect\ref{armineq}}{Proof of Proposition 3.3}}\label{secpfarm}

A significant step in the arguments of~\cite{Kesten87} is called the
``separation theorem'' (see also~\cite{Nolin}, Theorem 11). This states
roughly that, conditional
on the arm event $\Arm_{\si}(N,n)$, there is
probability bounded away from $0$ that arms
with the required colours can be found whose endpoints on the exterior
boundary of
the annulus are separated from one another by a given distance or more.
A formal statement
of the separation theorem is included in Section~\ref{secsep}
as Theorem~\ref{separation}; the proof is rather technical and very
similar to those
of~\cite{Kesten87,Nolin} and is therefore omitted.
The proof of Proposition~\ref{armineq} relies on arm-separation techniques.
More specifically, it
relies on Proposition~\ref{armextension}, which is an application of
the separation theorem,
Theorem~\ref{separation}.

The proof of Proposition~\ref{armineq} uses the following lemma, in
which the probability-vector
$\bp$ helps define the star--triangle transformations comprising the
map $\tau$.
%
%le3.4 #&#
\begin{lemma}\label{armtransfo}
Let $\Lattice=(V,E)$ be a mixed lattice,
let $\bp\in[0,1)^3$ be self-dual,
and let $k \in\{1, 2, 4, 6, \ldots\}$.
For $n/\sqrt3 > k +2$ and $\tau\in\{S^+ \comp T^+, S^- \comp T^-\}$,
we have (``surely'') that $\tau\Arm_k^{\Lattice}(n)
\subseteq\Arm_k^{\tau\Lattice}(n - 1)$.
\end{lemma}
\begin{pf}%[Proof of Lemma~\ref{armtransfo}]
Let $k \in\{1, 4, 6, \ldots\}$, we shall consider the case $k=2$ separately.
Let $\tau\in\{S^+ \comp T^+, S^- \comp T^-\}$
and $\omega\in\Arm_{k}^\Lattice(n)$.
Note that the points $x_r$, $r=0,1,\ldots\,$, are invariant under
$\tau$.

It is explained in Section~\ref{sectransformation} (see also
\cite{GM1}, Section 2)
that the image $\tau(\pi)$ of an $\omega$-open path $\pi$ contains
a $\tau(\omega)$-open
path of $\tau\Lat$ lying within distance $1$ of $\pi$.
Therefore, for $n/\sqrt3 > 2r+2$, if $C_{x_r}(\omega) \cap\pd B_n
\ne\es$,
then $C_{x_r}(\tau(\omega)) \cap\pd B_{n-1} \ne\es$.
The proof when $k=1$ is complete, and we\vspace*{1pt} assume now
that $k \ge4$. Let $j=k/2$ and $n/\sqrt3 > k+2$.
By Proposition~\ref{proppaths}, $x_r \nxlra{B_{n-1},\tau(\omega)}
x_s$ for $0 \le r < s\le j-1$,
whence $\tau(\omega)\in\Arm_k^{\tau\Lat}(n-1)$.

Finally, let $k=2$.
Let $\tau\in\{S^+ \comp T^+, S^- \comp T^-\}$
and $\omega\in\Arm_{2}^\Lattice(n)$.
Let $\Gamma$ (resp.,~$\Gamma^*$) be an open primal (resp., open$^*$
dual) path
starting at $x_0$ (resp., $y_0$)
that intersects~$\pd B_n$.
Since $x_0$ and $y_0$ are unchanged under $\tau$,
they are contained, respectively, in $\tau(\Gamma)$ and $\tau(\Gamma^*)$.
By the remarks in Section~\ref{sectransformation} concerning the
operation of $\tau$ on
open$^*$ dual paths, we conclude that $C_{x_0} \cap\pd B_{n-1} \ne\es
$ in $\tau\Lat$,
and similarly $C_{y_0}^* \cap\pd B_{n-1} \ne\es$ in $\tau\Lat^*$.
The proof is complete.
\end{pf}
\begin{pf*}{Proof of Proposition~\ref{armineq}}
Let $c$ and $N_1$ be as in Proposition~\ref{armextension}.
By making $n$ applications of\vadjust{\goodbreak} $\tau= S^+ \comp T^+$ to $\Lattice^0$,
we deduce that $\tau^n A_k^{\Lat^0}(2n) \subseteq A_k^{\Lat^n}(n)$.
Therefore, for $n \ge N_1$,
\begin{eqnarray*}
\PP_{(p_0, 1-p_0)}^\square[\Arm_k(n)]
& = & \PP_{\bp}^n [\Arm_{k}(n)] \\
&\geq& \PP_{\bp}^0 [\Arm_{k}(2n)]\qquad \mbox{by Lemma \ref
{armtransfo}} \\
& = & \PP_{\bp}^\triangle[\Arm_{k}(2n)] \\
& \geq & c \PP_{\bp}^\triangle[\Arm_{k}(n)] \qquad\mbox{by
Proposition~\ref{armextension}}.
\end{eqnarray*}
This proves the first inequality of Proposition~\ref{armineq}.

Let $n \geq\max\{k\sqrt3,N_1\}$, and consider the event $\Arm_k(n)$
on the lattice $\Lattice^n$.
If we apply $n$ times the transformation $S^- \comp T^-$ to $\Lattice^n$,
we obtain via Lemma~\ref{armtransfo} applied to the event $\Arm_{k}(2n)$
that:
\begin{eqnarray*}
\PP_{(p_0, 1-p_0)}^\square[\Arm_{k}(n)]
& = &\PP_{\bp}^n [\Arm_{k}(n)] \\
& \leq & c^{-1} \PP_{\bp}^n [\Arm_{k}(2n)] \qquad\mbox{by
Proposition~\ref{armextension}}\\
& \leq & c^{-1} \PP_{\bp}^0 [\Arm_{k}(n)] \qquad\mbox{by Lemma \ref
{armtransfo}}\\
& = & c^{-1} \PP_\bp^\triangle[\Arm_{k}(n)].
\end{eqnarray*}
The proof is complete.
\end{pf*}

%s3.5 #&#
\subsection{Separation theorem}\label{secsep}

The so-called ``separation theorem'' is a basic element in Kesten's work
on scaling relations in two dimensions.
It asserts roughly that, conditional on the occurrence of a given arm event,
there is probability bounded from $0$ that such arms may be found whose
endpoints
on the interior and exterior boundaries of the annulus are distant from
one another.
The separation theorem is useful since it permits the extensions of the arms
using box-crossings.

Kesten proved his theorem in~\cite{Kesten87} for homogeneous site
percolation models, while noting that
it is valid more generally. The proof has been reworked in
\cite{Nolin}, also in the context
of site percolation. The principal tool is the box-crossing property of
the critical
model. In this section, we state
a general form of the separation theorem, for use in both the current
paper and the forthcoming~\cite{GM3}.
The proof follows closely that found in
\cite{Kesten87,Nolin} and is omitted.

Let $G =(V,E)$ be a connected planar graph,
embedded in the plane in such a way that each edge is a straight
line segment, and let $\PP$ be a product measure on $\Omega= \{0,1\}^E$.
As usual we denote by $G^*$ the dual graph of $G$,
and more generally the superscript $*$
indicates quantities defined on the dual.
We shall use the usual notation from percolation
theory,~\cite{GrimmettPercolation},
and
we assume there exists a uniform upper bound $L<\oo$ on the lengths of
edges of $G$ and~$G^*$,
viewed as straight line segments of $\RR^2$.

The hypothesis required for the separation theorem concerns a lower
bound on the probabilities
of open and open$^*$ box-crossings.
Let $\omega\in\Om$ and let $R$ be a (nonsquare) rectangle of $\RR^2$.
A graph-path $\pi$ is said to \textit{cross} $R$
if $\pi$ contains an arc (termed a \textit{box-crossing}) lying in the
interior of $R$ except for its two endpoints, which are required
to lie, respectively, on the two shorter sides of $R$.
Note that box-crossings lie in the longer direction.
The rectangle $R$ is said to \textit{possess an open crossing}
(resp., \textit{open$^*$ dual crossing})
if there exists an open
path of $G$ (resp., open$^*$ path of $G^*$) crossing $R$,
and we write
$C(R)$ [resp., $C^*( R)$] for the event that this occurs.

Let $\sT$ be the set of translations of $\RR^2$, and $\tau\in\sT$.
Let $H_n=[0,2n]\times[0,n]$ and $V_n=[0,n]\times[0,2n]$, and
let $n_0=n_0(G)<\oo$ be minimal with the property that,
for all $\tau$ and all $n \ge n_0$, $\tau H_n$ and $\tau V_n $
possess crossings in both $G$ and $G^*$.
Let
%
%e3.4 #&#
%e3.5 #&#
\begin{eqnarray}
\label{G13}
b(G,\PP) &=& \inf\{ \PP(C(\tau H_n)),\PP(C( \tau V_n))\dvtx n \ge
n_0, \tau\in\sT\},\\
\label{G14}
b^*(G,\PP) &=& \inf\{ \PP(C^*(\tau H_n)),\PP(C^*(\tau V_n))\dvtx n
\ge n_0, \tau\in\sT\}
\end{eqnarray}
and
%
%e3.6 #&#
\begin{equation}\label{G1}
\beta= \beta(G,\PP) = \min\{b, b^*\}.
\end{equation}

The pair $(G, \PP)$ [resp., $(G^*,\PP^*)$] is said to have the
\textit
{box-crossing property}
if and only if $b(G,\PP)>0$ [resp.,
$b^*(G,\PP) > 0$].
In~\cite{GM1}, the box-crossing property is given in terms of boxes of
arbitrary
orientation and aspect-ratio. It is shown there that it suffices to
consider only horizontal and vertical
boxes. It is a consequence of the FKG inequality that
the box-crossing property does not depend on the chosen aspect-ratio
(so long as it is
strictly greater than $1$).
By~\cite{GM1}, Theorem 1.3, $\PP_\bp^\tri$ has the box-crossing
property whenever $\bp
\in[0,1)^3$
is self-dual.

Let $k \in\NN$ and $\si\in\{0,1\}^k$.
Rather than working with the arm events
$\Arm_{\sigma}(N,n)$ of Section~\ref{secexp}, we use instead the
events $\oArm_\si(N,n)$
defined in the same way except that $\La_n$ is replaced throughout the
definition by
$B_n=[-n,n]^2$, and arms are required to comprise edges that intersect $B_n$.
All constants in the following statements are permitted to depend on
the colour sequence $\sigma$.

Let $\oAnn(N,n)=B_n \setminus(-N,N)^2$ be the annulus with interior
boundary $\pd B_N$ and exterior boundary $\pd B_n$.
We shall consider open and open$^*$ crossings between
the interior and exterior boundaries.
We emphasize that the endpoints of these crossings are not required to
be graph vertices.

%f8 #&#
\begin{figure}

\includegraphics{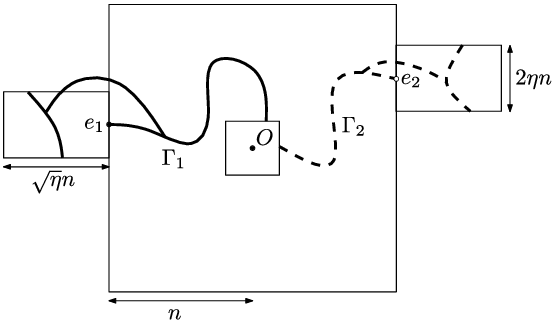}

\caption{A primal $\eta$-exterior-fence $\Ga_1$ with exterior
endpoint $e_1$,
and a dual $\eta$-exterior-fence $\Ga_2$.}
\label{figfence}
\end{figure}

For clarity, we concentrate first on the behaviour of crossings at
their exterior endpoints.
Let $\eta\in(0,1)$. A primal (resp., dual) \textit{$\eta$-exterior-fence}
is a set $\Ga$ of connected open (resp., open$^*$) paths comprising
the union of:
\begin{longlist}[(iii)]
\item[(i)] a crossing of $\oAnn(N,n)$ from its interior to its exterior
boundary, with exterior endpoint denoted
$\ext(\Ga)$,
\end{longlist}
together with certain further paths which we describe thus under the
assumption that
$\ext(\Ga)=(n,y)$ is on the right-hand side of $\pd B_n$:
\begin{longlist}[(iii)]
\item[(ii)] a vertical crossing of the box
$[n, (1+\sqrt{\eta})n] \times[y - \eta n, y + \eta n]$,
\item[(iii)] a connection between the above two crossings,
contained in $\ext(\Ga) +B_{\sqrt{\eta}n}$.
\end{longlist}
If $\ext(\Ga)$ is on a different side of $\pd B_n$, the event of
condition (ii) is replaced by
an appropriately rotated and translated event.
This definition is illustrated in Figure~\ref{figfence}.

One may similarly define an $\eta$-\textit{interior}-fence by
considering the
behaviour of the crossing near its interior endpoint.
We introduce also the concept of a primal (resp., dual)
$\eta$-\textit{fence}; this is a union of
an open (resp., open$^*$) crossing
of $\oAnn(N,n)$ together with further paths in the vicinities
of both interior and exterior endpoints along the lines of the above
definitions.

An \textit{$\eta$-landing-sequence} is a sequence of closed sub-intervals
$I = (I_i\dvtx i =1,2,\ldots, k)$ of $\pd B_1$,
taken in anticlockwise order,
such that each $I_i$ has length $\eta$, and
the minimal distance between any two intervals, and
between any interval and a corner of $B_1$,
is greater than $2\sqrt\eta$.
We shall assume that
%
%e3.7 #&#
\begin{equation}\label{eta}
0 < k\bigl(\eta+2\sqrt\eta\bigr) < 8,
\end{equation}
so that $\eta$-landing-sequences exist.

\begin{figure}

\includegraphics{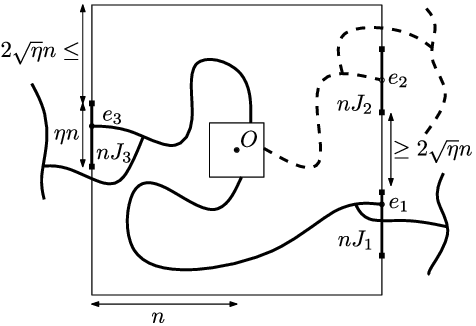}

\caption{The event $\Arm_{\sigma}^{\es,J} (N,n)$ with
$\sigma= (1,0,1)$ and $\eta$-landing-sequence $J$.
Each crossing $\Gamma_i$ is an $\eta$-exterior-fence
with exterior endpoint $e_i \in n J_i$.}
\label{figseparatedfences}
\end{figure}

Let $\eta,\eta'$ satisfy (\ref{eta}), and let $I$ (resp., $J$) be
an $\eta$-landing-sequence
(resp., $\eta'$-landing-sequence). Write $\oArm_{\sigma}^{I,J} (N,n)$
for the event that there exists a sequence of
$\eta$-fences $(\Gamma_i\dvtx i =1,2, \ldots, k)$
in the annulus $B_n \setminus(-N,N)^2$,
with colours prescribed by $\sigma$, such that, for all $i$,
the interior (resp., exterior) endpoint of $\Ga_i$ lies in $N I_i$
(resp., $n J_i$).
Let $\oArm_{\sigma}^{I,\es} (N,n)$ [resp., $\oArm_{\sigma}^{\es
,J} (N,n)$]
be given similarly in terms of
$\eta$-interior-fences (resp., $\eta'$-exterior-fences).
Note that
%
%e3.8 #&#
\begin{equation}\label{will1}
\oArm_{\sigma}^{I,J} (N,n) \subseteq\oArm_\si^{\es,J}(N,n),\qquad \oArm
_\si^{I,\es}(N,n) \subseteq\oArm_\si(N,n).
\end{equation}
These definitions are illustrated in Figure~\ref{figseparatedfences}.

%f9 #&#

In the proof of the forthcoming Proposition~\ref{armextension}
(and nowhere else), we shall make use of
a piece of related notation introduced here. Let $k=2j \ge2$,
and let $\eta$ and $I$ be as above. As explained
in the proof of Lemma~\ref{lemhook}, the event $\Arm_k(n)$ of
Section~\ref{secmainp} requires the existence of
an alternating sequence of open and open$^*$ paths joining the set $\{
x_0,y_0,x_1,y_1,\ldots,x_{j-1}\}$
to the boundary $\pd B_n$.
Let $\Arm_k^{I} (n)$ be the sub-event
in which the exterior endpoints of the open (resp., open$^*$)
paths lie in $I_1,I_3,\ldots,I_{k-1}$ (resp., $I_2,I_4,\ldots,I_k$),
and in addition these
exterior endpoints have associated paths as given in
(ii) and (iii) of the above definition of an $\eta$-exterior-fence.

We now state the separation theorem.
The proof is omitted, and may be
constructed via careful readings
of the appropriate sections of~\cite{Kesten87,Nolin}.
There is a small complication arising from the fact that the endpoints
of box-crossings
are not necessarily vertices of the relevant graph, and this is
controlled using the uniform
upper bound $L$ on the lengths of embeddings of edges.
%
%th3.5 #&#
\begin{theorem}[(Separation theorem)]\label{separation}
Let $k \in\NN$, and $\si\in\{0,1\}^k$.
For $\beta_0 >0$, $M\in\NN$, and $\eta_0>0$,
there exist constants $c > 0$ and $n_1 \in\NN$ such that:
for any pair $(G,\PP)$ with $\beta(G,\PP) > \beta_0$ and $n_0(G)\le M$,
for all $\eta,\eta' > \eta_0$ satisfying (\ref{eta}),
all $\eta$-landing-sequences $I$ and $\eta'$-landing-sequences $J$,
and all $N \geq n_1$ and $n \geq2N$, we have
\[
\PP[ \oArm_{\sigma}^{I,J} (N,n) ] \geq
c \PP[ \oArm_{\sigma} (N,n) ].
\]
\end{theorem}

Amongst the consequences of Theorem~\ref{separation} is the following.
The proof
(also omitted) is essentially that of~\cite{Nolin}, Proposition 12,
and it uses the extension of paths by judiciously positioned
box-crossings.
%
%co3.6 #&#
\begin{cor}\label{corsepcor}
Suppose that $\beta= \beta(G,\PP) > 0$.
For $k \in\NN$ and $\sigma\in\{0,1\}^k$,
there exists $c =c(\beta,\si) > 0$ and $N_0 \in\NN$ such that, for all
$N \ge N_0$ and $n\ge2N$,
\[
\PP[\oArm_{\sigma}(N,2n)] \geq c \PP[\oArm_{\sigma}(N,n)].
\]
\end{cor}

The proof of Proposition~\ref{armineq} makes use of the following application
of Corollary~\ref{corsepcor} to the pairs $(\Lat^m, \PP_\bp^m)$.
%
%pr3.7 #&#
\begin{prop}\label{armextension}
For $k \in\{1, 2, 4, 6, \ldots\}$ and
a self-dual triplet $\bp\in[0,1)^3$ with $p_0>0$,
there exist $c > 0$ and $N_1 \in\NN$ such that, for $m \ge0$
and $n\ge N_1$,
\[
\PP_\bp^m [\Arm_{k}(2n)] \geq c \PP_\bp^m [\Arm_{k}(n)].
\]
\end{prop}
\begin{pf}
The box-crossing property has been studied in~\cite{GM1} in the context
of hybrids of
$\ZZ^2$/$\TT$ or $\ZZ^2$/$\HH$ type. The arguments of~\cite{GM1}
may be adapted as follows
to obtain that, for given self-dual $\bp\in[0,1)^3$ with $p_0>0$,
the pairs
$(\Lat^m,\PP_\bp^m)$ satisfy a \textit{uniform} box-crossing
property in
the sense that:
there exists $\beta_0>0$ such that
%
%e3.9 #&#
\begin{equation}\label{G15}
\beta(\Lat^m,\PP_\bp^m) > \beta_0,\qquad m \ge0.
\end{equation}

Write $B_{M,N} = [0,M] \times[0,N]$, and denote by $\Ch(B_{M,N})$
[resp.,
$\Cv(B_{M,N})$] the event that there exists a horizontal
(resp., vertical) open crossing of $B_{M,N}$ (with a similar
notation $\Ch^*$, $\Cv^*$ for dual crossings).
Since
every translate of $B_{M,3N}$ contains a rectangle with dimensions
$M\times N$
lying in either the square or triangular part of $\Lat^m$,
\[
\PP_\bp^m[\Ch(\tau B_{M,3N}) ]\ge\min\bigl\{\PP_\bp^\tri[ \Ch
(B_{M,N})],
\PP_{(p_0, 1-p_0)}^\squ[ \Ch(B_{M,N})]\bigr\}
\]
for all $\tau\in\sT$.
The dual model lives on a mixed square/hexagonal lattice with parameter
$1-\bp$, and the same
inequality holds with $\Ch$ replaced by $\Ch^*$.
By~\cite{GM1}, Theorem 1.3, there exists $b_1=b_1(\bp)>0$ such that
%
%e3.10 #&#
\begin{equation}\label{G10}
\PP_\bp^m[\Ch(\tau B_{M,3N})], \PP_\bp^m[\Ch^*(\tau B_{M,3N})]
\ge b_1,\qquad m\ge0, M,N\ge1, \tau\in\sT.\hspace*{-38pt}
\end{equation}

Adapting the proof of~\cite{GM1}, Proposition 3.8, we obtain that
\[
\PP^m_\bp[ \Cv(\tau B_{3N,N}) ] \geq
\PP_\bp^\tri[ \Cv(B_{N,2N}) ],\qquad m\ge0, N \ge1, \tau\in
\sT.
\]
The same inequality holds with $\Cv$ replaced by $\Cv^*$, as above,
and therefore there
exists $b_2=b_2(\bp)>0$ such that
%
%e3.11 #&#
\begin{equation}\label{G11}\qquad
\PP_\bp^m[\Cv(\tau B_{3N,N})], \PP_\bp^m[\Cv^*(\tau B_{3N,N})]
\ge b_2,\qquad m\ge0, N\ge1, \tau\in\sT.
\end{equation}
Inequalities (\ref{G10}) and (\ref{G11}) imply as in the proof of
\cite{GM1}, Proposition 3.1, that (\ref{G15}) holds for some $n_0<\oo$ and
$\beta_0>0$,
and we choose these accordingly.

Let $\si$ be an alternating colour sequence of length $k$ (we set $\si
=\{1\}$
when $k=1$), and denote $\oArm_\si(N,n)$ by $\oArm_k(N,n)$.

Let $\eta$ satisfy (\ref{eta}) and let $I$ be an $\eta$-landing sequence.
Let $c=c_2$ and $n_1$ be as in Theorem~\ref{separation}.
By Corollary~\ref{corsepcor}, there exists $c_0=c_0(\beta_0,k)>0$
and $N_0 \ge n_1$ such that
%
%e3.12 #&#
\begin{equation}\label{G16}
\PP_\bp^m [\oArm_k(N,2n)] \geq c_0 \PP_\bp^m [\oArm_k(N,n)],\qquad
m \ge0, n \ge2N \ge2N_0.
\end{equation}
Therefore,
%
%e3.13 #&#
\begin{equation}\label{G17}
\PP_\bp^m[\Arm_k(n)] \le\PP_\bp^m[\oArm_k(N,n)] \le
c_0^{-1}\PP_\bp^m[\oArm_k(N,2n)].
\end{equation}

By an elementary consideration of
paths of $\Lat^m$, there exist $N_1\ge2 N_0$ and $c_1=c_1(\bp,N_1)$
such that
%
%e3.14 #&#
\begin{equation}
\label{G18}
\PP_\bp^m \bigl[\Arm_k^{I}\bigl(\tfrac12 N_1\bigr)\bigr] \ge c_1,\qquad m \ge0.
\end{equation}
By Theorem~\ref{separation} and (\ref{will1}),
%
%e3.15 #&#
\begin{equation}\label{G19}
\PP_\bp^m[\oArm_k^{I,\es}(N_1,2n)] \ge c_2
\PP_\bp^m[\oArm_k(N_1,2n)],\qquad
m \ge0, n \ge N_1.
\end{equation}
Furthermore, by the uniform box-crossing property as in~\cite{Nolin},
Proposition~12, there exists $c_3=c_3(\bp,\eta)>0$ such that
\begin{eqnarray*}
\PP_\bp^m[\Arm_k(2n)] &\ge& c_3\PP_\bp^m\bigl[\Arm_k^{I}\bigl(\tfrac12
N_1\bigr)\bigr] \PP_\bp^m[\oArm_k^{I,\es}(N_1,2n)]\\
&\ge& c_1c_2c_3 \PP_\bp^m[\oArm_k(N_1,2n)] \qquad\mbox{by (\ref
{G18}) and (\ref{G19})}.
\end{eqnarray*}
The claim follows by (\ref{G17}) with $c=c_0c_1c_2c_3$.
\end{pf}

\section*{Acknowledgment}

Alexander Holroyd proposed a consideration of the Archimedean
lattices.

%suskaldyti doi

% imsref loaded by lrinkeviciute, 2012-05-30 13:50:48
% imsref loaded by lrinkeviciute, 2012-05-30 13:57:56

\printaddresses


\begin{thebibliography}{21}
% BibTex style file: ims.bst, 2011-05-30
% Default style options (sort=0,type=number).
% Used options (sort=1,type=number).

%b1 ###
\bibitem{ADA}
\begin{barticle}[auto:STB|2012/05/28|15:16:20]
\bauthor{\bsnm{Aizenman},~\bfnm{M.}\binits{M.}},
  \bauthor{\bsnm{Duplantier},~\bfnm{B.}\binits{B.}} \AND
  \bauthor{\bsnm{Aharony},~\bfnm{A.}\binits{A.}}
(\byear{1999}).
\btitle{Path-crossing exponents and the external perimeter in 2{D}
  percolation}.
\bjournal{Phys. Rev. Lett.}
\bvolume{83}
\bpages{1359--1362}.
\bptok{imsref}%
\end{barticle}
\endbibitem

%b2 ###
\bibitem{BR}
\begin{bincollection}[mr]
\bauthor{\bsnm{Bollob{\'a}s},~\bfnm{B{\'e}la}\binits{B.}} \AND
  \bauthor{\bsnm{Riordan},~\bfnm{Oliver}\binits{O.}}
(\byear{2010}).
\btitle{Percolation on self-dual polygon configurations}.
In \bbooktitle{An Irregular Mind}.
\bseries{Bolyai Society Mathematical Studies}
\bvolume{21}
\bpages{131--217}.
\bpublisher{J\'anos Bolyai Math. Soc.}, \baddress{Budapest}.
\bid{doi={10.1007/978-3-642-14444-8_3}, mr={2815602}}
\bptok{imsref}%
\end{bincollection}
\endbibitem

%b3 ###
\bibitem{BdTB}
\begin{bmisc}[auto:STB|2012/05/28|15:16:20]
\bauthor{\bsnm{Boutillier},~\bfnm{C.}\binits{C.}} \AND \bauthor{\bparticle{de}
  \bsnm{Tili{\`e}re},~\bfnm{B.}\binits{B.}}
(\byear{2010}).
\bhowpublished{Statistical mechanics on isoradial graphs. Available at
  arXiv:\arxivurl{1012.2955}.}
\bptok{imsref}%
\end{bmisc}
\endbibitem

%b4 ###
\bibitem{GrimmettPercolation}
\begin{bbook}[mr]
\bauthor{\bsnm{Grimmett},~\bfnm{Geoffrey}\binits{G.}}
(\byear{1999}).
\btitle{Percolation},
\bedition{2nd} ed.
\bseries{Grundlehren der Mathematischen Wissenschaften [Fundamental Principles
  of Mathematical Sciences]}
\bvolume{321}.
\bpublisher{Springer}, \baddress{Berlin}.
\bid{mr={1707339}}
\bptok{imsref}%
\end{bbook}
\endbibitem

%b5 ###
\bibitem{GM3}
\begin{bmisc}[auto:STB|2012/05/28|15:16:20]
\bauthor{\bsnm{Grimmett},~\bfnm{G.~R.}\binits{G.~R.}} \AND
  \bauthor{\bsnm{Manolescu},~\bfnm{I.}\binits{I.}}
(\byear{2011}).
\bhowpublished{Bond percolation on isoradial graphs.
Preprint. Available at \arxivurl{arXiv:1204.0505}.}
\bptok{imsref}%
\end{bmisc}
\endbibitem

%b6 ###
\bibitem{GM1}
\begin{barticle}[auto:STB|2012/05/28|15:16:20]
\bauthor{\bsnm{Grimmett},~\bfnm{G.~R.}\binits{G.~R.}} \AND
  \bauthor{\bsnm{Manolescu},~\bfnm{I.}\binits{I.}}
(\byear{2013}).
\btitle{Inhomogeneous bond percolation on the square, triangular, and
  hexagonal lattices}.
  \bjournal{Ann. Probab.}
  \bvolume{41}
  \bpages{2990--3025}.
\bptok{imsref}%
\end{barticle}
\endbibitem

%b7 ###
\bibitem{GS}
\begin{bbook}[mr]
\bauthor{\bsnm{Gr{\"u}nbaum},~\bfnm{Branko}\binits{B.}} \AND
  \bauthor{\bsnm{Shephard},~\bfnm{G.~C.}\binits{G.~C.}}
(\byear{1987}).
\btitle{Tilings and Patterns}.
\bpublisher{Freeman}, \baddress{New York}.
\bid{mr={0857454}}
\bptok{imsref}%
\end{bbook}
\endbibitem

%b8 ###
\bibitem{Ken02}
\begin{bincollection}[mr]
\bauthor{\bsnm{Kenyon},~\bfnm{Richard}\binits{R.}}
(\byear{2004}).
\btitle{An introduction to the dimer model}.
In \bbooktitle{School and {C}onference on {P}robability {T}heory}.
\bseries{ICTP Lect. Notes}
\bvolume{XVII}
\bpages{267--304 (electronic)}.
\bpublisher{Abdus Salam Int. Cent. Theoret. Phys.}, \baddress{Trieste}.
\bid{mr={2198850}}
\bptok{imsref}%
\end{bincollection}
\endbibitem

%b9 ###
\bibitem{KenS}
\begin{barticle}[mr]
\bauthor{\bsnm{Kenyon},~\bfnm{Richard}\binits{R.}} \AND
  \bauthor{\bsnm{Schlenker},~\bfnm{Jean-Marc}\binits{J.-M.}}
(\byear{2005}).
\btitle{Rhombic embeddings of planar quad-graphs}.
\bjournal{Trans. Amer. Math. Soc.}
\bvolume{357}
\bpages{3443--3458 (electronic)}.
\bid{doi={10.1090/S0002-9947-04-03545-7}, issn={0002-9947}, mr={2146632}}
\bptok{imsref}%
\end{barticle}
\endbibitem

%b10 ###
\bibitem{Kestenbook}
\begin{bbook}[mr]
\bauthor{\bsnm{Kesten},~\bfnm{Harry}\binits{H.}}
(\byear{1982}).
\btitle{Percolation Theory for Mathematicians}.
\bseries{Progress in Probability and Statistics}
\bvolume{2}.
\bpublisher{Birkh\"auser}, \baddress{Boston}.
\bid{mr={0692943}}
\bptok{imsref}%
\end{bbook}
\endbibitem

%b11 ###
\bibitem{Kes87a}
\begin{bincollection}[mr]
\bauthor{\bsnm{Kesten},~\bfnm{Harry}\binits{H.}}
(\byear{1987}).
\btitle{A scaling relation at criticality for {$2$}{D}-percolation}.
In \bbooktitle{Percolation Theory and Ergodic Theory of Infinite Particle
  Systems ({M}inneapolis, {M}inn., 1984--1985)}.
\bseries{The IMA Volumes in Mathematics and its Applications}
\bvolume{8}
\bpages{203--212}.
\bpublisher{Springer}, \baddress{New York}.
\bid{doi={10.1007/978-1-4613-8734-3_12}, mr={0894549}}
\bptok{imsref}%
\end{bincollection}
\endbibitem

%b12 ###
\bibitem{Kesten87}
\begin{barticle}[mr]
\bauthor{\bsnm{Kesten},~\bfnm{Harry}\binits{H.}}
(\byear{1987}).
\btitle{Scaling relations for {$2$}{D}-percolation}.
\bjournal{Comm. Math. Phys.}
\bvolume{109}
\bpages{109--156}.
\bid{issn={0010-3616}, mr={0879034}}
\bptok{imsref}%
\end{barticle}
\endbibitem

%b13 ###
\bibitem{Nolin}
\begin{barticle}[mr]
\bauthor{\bsnm{Nolin},~\bfnm{Pierre}\binits{P.}}
(\byear{2008}).
\btitle{Near-critical percolation in two dimensions}.
\bjournal{Electron. J. Probab.}
\bvolume{13}
\bpages{1562--1623}.
\bid{doi={10.1214/EJP.v13-565}, issn={1083-6489}, mr={2438816}}
\bptok{imsref}%
\end{barticle}
\endbibitem

%b14 ###
\bibitem{PW}
\begin{barticle}[mr]
\bauthor{\bsnm{Parviainen},~\bfnm{Robert}\binits{R.}} \AND
  \bauthor{\bsnm{Wierman},~\bfnm{John~C.}\binits{J.~C.}}
(\byear{2005}).
\btitle{Inclusions and non-inclusions among the {A}rchimedean and {L}aves
  lattices, with applications to bond percolation thresholds}.
\bjournal{Congr. Numer.}
\bvolume{176}
\bpages{89--128}.
%  Theory, and Computing}.
\bid{issn={0384-9864}, mr={2198637}}
\bptok{imsref}%
\end{barticle}
\endbibitem

%b15 ###
\bibitem{SedW}
\begin{barticle}[mr]
\bauthor{\bsnm{Sedlock},~\bfnm{Matthew R.~A.}\binits{M.~R.~A.}} \AND
  \bauthor{\bsnm{Wierman},~\bfnm{John~C.}\binits{J.~C.}}
(\byear{2009}).
\btitle{Equality of bond-percolation critical exponents for pairs of dual
  lattices}.
\bjournal{Phys. Rev. E (3)}
\bvolume{79}
\bpages{051119, 10}.
\bid{doi={10.1103/PhysRevE.79.051119}, issn={1539-3755}, mr={2551408}}
\bptok{imsref}%
\end{barticle}
\endbibitem

%b16 ###
\bibitem{Smirnov}
\begin{barticle}[mr]
\bauthor{\bsnm{Smirnov},~\bfnm{Stanislav}\binits{S.}}
(\byear{2001}).
\btitle{Critical percolation in the plane: Conformal invariance, {C}ardy's
  formula, scaling limits}.
\bjournal{C. R. Acad. Sci. Paris S\'er. I Math.}
\bvolume{333}
\bpages{239--244}.
\bid{doi={10.1016/S0764-4442(01)01991-7}, issn={0764-4442}, mr={1851632}}
\bptok{imsref}%
\end{barticle}
\endbibitem

%b17 ###
\bibitem{Smirnov-Werner}
\begin{barticle}[mr]
\bauthor{\bsnm{Smirnov},~\bfnm{Stanislav}\binits{S.}} \AND
  \bauthor{\bsnm{Werner},~\bfnm{Wendelin}\binits{W.}}
(\byear{2001}).
\btitle{Critical exponents for two-dimensional percolation}.
\bjournal{Math. Res. Lett.}
\bvolume{8}
\bpages{729--744}.
\bid{issn={1073-2780}, mr={1879816}}
\bptok{imsref}%
\end{barticle}
\endbibitem

%b18 ###
\bibitem{WernerSMF}
\begin{bbook}[auto:STB|2012/05/28|15:16:20]
\bauthor{\bsnm{Werner},~\bfnm{W.}\binits{W.}}
(\byear{2009}).
\btitle{Percolation et Mod\`ele D'Ising}.
\bseries{Cours Specialis\'es}
\bvolume{16}.
\bpublisher{Soci\'et\'e Math\'ematique de France}, \baddress{Paris}.
\bptok{imsref}%
\end{bbook}
\endbibitem

%b19 ###
\bibitem{Z1}
\begin{barticle}[mr]
\bauthor{\bsnm{Ziff},~\bfnm{Robert~M.}\binits{R.~M.}}
(\byear{2006}).
\btitle{Generalized cell--dual-cell transformation and exact thresholds for
  percolation}.
\bjournal{Phys. Rev. E (3)}
\bvolume{73}
\bpages{016134, 6}.
\bid{doi={10.1103/PhysRevE.73.016134}, issn={1539-3755}, mr={2223061}}
\bptok{imsref}%
\end{barticle}
\endbibitem

%b20 ###
\bibitem{ZScull}
\begin{barticle}[mr]
\bauthor{\bsnm{Ziff},~\bfnm{Robert~M.}\binits{R.~M.}} \AND
  \bauthor{\bsnm{Scullard},~\bfnm{Christian~R.}\binits{C.~R.}}
(\byear{2006}).
\btitle{Exact bond percolation thresholds in two dimensions}.
\bjournal{J. Phys. A}
\bvolume{39}
\bpages{15083--15090}.
\bid{doi={10.1088/0305-4470/39/49/003}, issn={0305-4470}, mr={2277091}}
\bptok{imsref}%
\end{barticle}
\endbibitem

\end{thebibliography}
\end{document}